\documentclass[11pt,a4paper]{amsart}

\usepackage[T1]{fontenc}
\usepackage{lmodern}
\usepackage[utf8]{inputenc}
\usepackage[margin=1.15in]{geometry}
\usepackage{mathtools,amssymb,amsthm,amsfonts,mathrsfs,bm, bbm}
\usepackage{microtype}
\usepackage[dvipsnames]{xcolor}
\usepackage{hyperref}
\usepackage[nameinlink,capitalise,noabbrev]{cleveref}

\usepackage[
  backend=biber,
  style=numeric,
  sorting=nyt,
  giveninits=true,
  maxbibnames=99
]{biblatex}

\usepackage{graphicx}
\usepackage{tikz-cd}
\usepackage{tikz-3dplot}
\tikzset{smalltext/.style={"\textup{\small #1}" description}}

\hypersetup{
  colorlinks=true,
  linkcolor=blue!50!black,
  citecolor=blue!50!black,
  urlcolor=blue!50!black
}

\numberwithin{equation}{section}

\theoremstyle{plain}
\newtheorem{theorem}{Theorem}[section]
\newtheorem{proposition}[theorem]{Proposition}
\newtheorem{lemma}[theorem]{Lemma}
\newtheorem{corollary}[theorem]{Corollary}

\theoremstyle{definition}

\theoremstyle{remark}
\newtheorem{remark}[theorem]{Remark}

\newcommand{\E}{\mathbb{E}}
\renewcommand{\H}{\mathcal{H}}
\newcommand{\N}{\mathbb{N}}
\newcommand{\R}{\mathbb{R}}
\newcommand{\T}{\mathbb{T}}
\newcommand{\Z}{\mathbb{Z}}

\newcommand{\eps}{\varepsilon}

\renewcommand{\P}{P}

\newcommand{\jb}[1]{\langle #1\rangle}

\newcommand{\EV}[1]{\mathbb{E}\bracket*{#1}}
\newcommand{\EVMU}[1]{\mathbb{E}_\mu \bracket*{#1}}

\newcommand{\wick}[1]{:\! #1 \!:}

\newcommand{\les}{\lesssim}
\DeclarePairedDelimiter{\bracket}{[}{]}

\title[Uniqueness of invariant measures]
  {Uniqueness of invariant measures for mean-field stochastic damped wave and heat equations}
\author{Alexander Czernik}
\address{Université Jean Monnet Saint-Étienne}
\email{alexander.czernik@univ-st-etienne.fr}
\date{}
\subjclass[2020]{60H15, 35L71, 37L40}
\keywords{stochastic damped wave equation, invariant measure,
  mean-field nonlinearity, Wick renormalization}

\begin{document}

\begin{abstract}
We study invariant measures for a Wick-renormalized mean-field stochastic damped wave equation on the two-dimensional torus and for the corresponding stochastic heat equation. With Wick ordering relative to the massive Gaussian free field with covariance \((1-\Delta)^{-1}\), we prove uniqueness among invariant laws whose renormalized second moment belongs to \(H^{-\sigma}(\mathbb T^2)\), where \(0<\sigma<\tfrac12\). Stationarity reduces the mean-field equation to a linear equation with a deterministic distribution-valued potential. We show that self-consistency forces this deterministic potential to vanish. The damped wave equation therefore has invariant law \(\mathcal N(0,(1-\Delta)^{-1})\otimes\mathcal N(0,\varepsilon^{-2}I)\), with covariance understood distributionally, while the heat equation has the Gaussian free field as its invariant law. This resolves self-consistency for spatially varying distribution-valued potentials and establishes uniqueness in the specified class.
As the inertia parameter tends to zero, the position marginal remains equal to the invariant law of the corresponding stochastic heat equation, whereas the velocity marginal assigns vanishing probability to every bounded subset of \(H^{-1-\sigma}(\mathbb T^2)\), establishing Smoluchowksi-Kramers type behavior. 
\end{abstract}

\maketitle

\section{Introduction}

Let $\xi$ be space-time white noise. We consider here the following stochastic damped wave equation with a wick-renormalized mean field nonlinearity (SDMFW) on the torus $\T^2 = \R^2 / (2\pi \Z)^2$ 
\begin{align}\tag{SDMFW}\label{eq:SDMFW}
\begin{cases}
\eps^2 \partial_t^2 u_\eps +  \partial_t u_\eps + (1 - \Delta) u_\eps + \E[\wick{u_\eps^2}]u_\eps   = \sqrt{2}\xi,    \\
(u_\eps,\partial_t u_\eps)_{| t=0} =(\phi _0,\phi _1) 
\end{cases}
\quad (t,x ) \in \R_+ \times \T^2
\end{align}
for $\eps \in (0,1]$. The Wick-renormalization is defined in Section~\ref{sec:SDMFW_inv_msr}.

We investigate the dynamics under the assumption of the existence of an invariant measure $\mu$ for the flow $u$. Restricting the class of measures to those who satisfy a specific moment condition, we show that such a measure is unique and show that this measure is the centered Gaussian measure with distributional covariance as specified in Section~\ref{sec:SDMFW_inv_msr}. This goal is summarized in Theorem~\ref{theorem:uniqueness_SDMFW}
\begin{theorem}\label{theorem:uniqueness_SDMFW}
    Given $\sigma \in (0,\frac{1}{2})$, $\eps \in (0,1]$ and let $\mathcal{H}^{-\sigma} \coloneq H^{-\sigma}(\T^2) \times H^{-1-\sigma}(\T^2)$. 
    Then $\rho_\eps = \mathcal{N}(0, \Lambda^{-1}) \otimes \mathcal{N}(0, \eps^{-2}I)$ is the unique invariant measure of the~\eqref{eq:SDMFW} in the class
    \[
        \mathcal W
        :=\left\{
            \mu\in\mathcal P(\mathcal H^{-\sigma}):
            \mathbb E_\mu[\wick{u^2}]\in H^{-\sigma}(\T^2)
        \right\}.
    \]
\end{theorem}

Furthermore, we generalize this uniqueness result to the parabolic stochastic mean-field heat equation (SMFH):
\begin{align}\tag{SMFH}\label{eq:SMFH}
\begin{cases}
 \partial_t u + (1 - \Delta) u + \E[\wick{u^2}]u   = \sqrt{2}\xi,    \\
u_{| t=0} =\phi _0
\end{cases}
\quad (t,x ) \in \R_+ \times \T^2.
\end{align}
Finally, we demonstrate a Smoluchowski-Kramers type behavior for the invariant measures $\rho_\eps$ and $\rho_0$ as $\eps \to 0$.
This result can be directly generalized to the unique invariant measure of the stochastic damped wave equation with cubic nonlinearity (see \cite{tolomeo2023ergodicityhyperbolicpphi2model}) and the invariant measure of the stochastic heat equation with cubic nonlinearity (see \cite{TsatsoulisWeber2018}). The reason is that the proof concerns effectively only the second marginal which is the same for both equation families. 
\begin{corollary}\label{theorem:convergence_uniqueness_SMFH}
    Given $\sigma \in (0,\frac{1}{2})$. Then $\rho_0 = \mathcal{N}(0, (1-\Delta)^{-1})$ is the unique invariant measure of the \eqref{eq:SMFH} in the class 
    \begin{align*}
        \mathcal{W}_0:=\left\{
            \mu\in\mathcal P(H^{-\sigma}(\T^2)):
            \mathbb E_\mu[\wick{u^2}]\in H^{-\sigma}(\T^2)
        \right\}.
    \end{align*}
    Furthermore, for $\sigma > 0$, 
    \begin{align*}
        \rho_\eps \xrightarrow[]{\mathrm{rv}} \rho_0 \quad \text{as $\eps \to 0$}
    \end{align*}
    in the sense that projected onto its first component $(\pi_1)_* \rho_\eps = \rho_0$ for all $\eps \in (0,1]$ and for all radial vanishing functions
    \begin{align*}
        f \in C_{\mathrm{rv}}(H^{-1-\sigma}(\T^2))
    \coloneqq
    \left\{
        f\in C_b(H^{-1-\sigma}(\T^2)) \mid
        \lim_{R\to\infty}
        \sup_{\|x\|_X\geq R}|f(x)|=0
    \right\}
    \end{align*}
    the second component is vanishing:
    \begin{align*}
        \int_{H^{-1-\sigma}(\T^2)}
        f(v)\,d((\pi_2)_*\rho_\eps)(v)
        \longrightarrow 0 \quad \text{as $\eps \to 0$}.
    \end{align*}
\end{corollary}

Our interest for the \eqref{eq:SMFH} and the \eqref{eq:SDMFW} arises from the stochastic quantization program proposed by Parisi and Wu in 1981 \cite{Parisi1981}, originally proposed as a method to evaluate path integrals in gauge theories without the need for gauge-fixing. The idea is that given a target Euclidean quantum field measure, we construct Langevin dynamics whose solution converges in Law to the target measure. Then, in principle, we can sample from the target measure via a Monte Carlo Method. 
Mathematically, in recent years the Euclidean $\Phi^4$ model and different approaches of its stochastic quantization have attracted a lot of attention. 
The most straight forward stochastic quantization is given by 
\begin{align}\label{eq:SNLH_unrenormalized}
\begin{cases}
 \partial_t u + (1 - \Delta) u + u^3  = \sqrt{2}\xi,    \\
u_{| t=0} =\phi _0.
\end{cases}
\end{align}
Since this is a parabolic equation we refer to this approach as the parabolic stochastic quantization of the $\Phi^4$ model. 
For $d = 2$, this constitutes a singular stochastic partial differential equation (SPDE), consequently, to make sense of this equation a renormalization is necessary. Considering the $\Phi^4_2$ Euclidean quantum field theory, the $\Phi^4$ model restricted to two dimensions, this normalization corresponds to the second moment of the free field theory, namely the second moment of the Gaussian measure $\mathcal{N}(0, (1-\Delta)^{-1})$ which formally equals $\infty$ (see \cite{GlimmJaffe1987}). 

The construction of the $\Phi^4_2$ model traces back to \cite{GlimmJaffe1987}. For the associated parabolic dynamics \cite{DaPratoDebussche2003} established global strong solutions on $\mathbb{T}^2$ for Gibbs-almost-every initial condition. \cite{MourratWeber2017Plane} subsequently proved global well-posedness for arbitrary initial conditions in suitable distribution spaces, both on the torus and on $\mathbb{R}^2$. On $\mathbb{T}^2$, unique ergodicity of the $\Phi^4_2$ measure was established in \cite{TsatsoulisWeber2018}, completing the parabolic stochastic quantization of the $\Phi^4_2$-model on the two-dimensional torus. 


An alternative approach to the parabolic stochastic quantization of $\Phi^4$ is provided by the hyperbolic stochastic quantization, given by the second-order equation
\begin{align}\label{eq:SNLDW_unrenormalized}
\begin{cases}
\partial_t^2 u+\partial_t u+(1-\Delta)u+u^3=\sqrt{2}\xi,\\
(u,\partial_t u)|_{t=0}=(\phi_0,\phi_1).
\end{cases}
\end{align}
Its Hamiltonian structure, with position variable $u$ and velocity variable $\partial_t u$, supplemented by damping and noise, connects this approach with sampling methods based on momentum variables, including Hamiltonian Monte Carlo \cite{DuaneEtAl1987}.

In two dimensions, Wick renormalization permits the construction of these dynamics despite the absence of parabolic smoothing.
On $\mathbb{T}^2$, \cite{GKOT2022} established global well-posedness for Gibbs-almost-every initial condition and invariance of the Gibbs measure. On $\mathbb{R}^2$, invariant Gibbs dynamics were constructed by \cite{OTWZ2026}. While on the torus, \cite{Tolomeo2023Ergodicity} proved ergodicity of the Gibbs measure and uniqueness within a class of invariant measures satisfying a moment condition. In particular, ergodicity yields almost-sure convergence of time averages for Gibbs-almost-every initial condition, partly completing the hyperbolic stochastic quantization program for the $\Phi^4_2$ model on the two-dimensional torus.

Equation \eqref{eq:SDMFW} arises as the mean-field limit of the hyperbolic stochastic quantization of the $O(N)$ linear sigma model on $\mathbb{T}^2$, an $N$-component generalization of the $\Phi^4_2$ model considered by \cite{Wilson1973}. The corresponding system consists of $N$ stochastic damped wave equations coupled through a Wick-renormalized cubic interaction with mean-field scaling $1/N$.

For the parabolic dynamics \cite{Shen2022} established convergence to the mean-field equation \eqref{eq:SMFH} in probability in suitable path spaces and under appropriate assumptions on the initial data. For sufficiently large mass, they also proved the Gaussian free field measure is the unique invariant measure of the \eqref{eq:SMFH} and convergence of every fixed collection of Gibbs marginals to the corresponding product of Gaussian free field measures. \cite{DelgadinoSmith2026} subsequently removed the large-mass restriction for this static convergence.

\cite{liu2025} established global well-posedness and the mean-field limit for the hyperbolic dynamics, again with convergence in probability on finite time intervals under suitable initial-data assumptions. They also proved convergence of the stationary Gibbs dynamics for every positive mass. At mass one, the limiting stationary law equals $\rho_1$.
This work does not establish uniqueness of the invariant law of the hyperbolic mean-field equations. 

In Theorem \ref{theorem:uniqueness_SDMFW} we establish uniqueness at mass one under the stated moment condition, complementing the large-mass uniqueness result for the parabolic mean-field dynamics in \cite{Shen2022}.

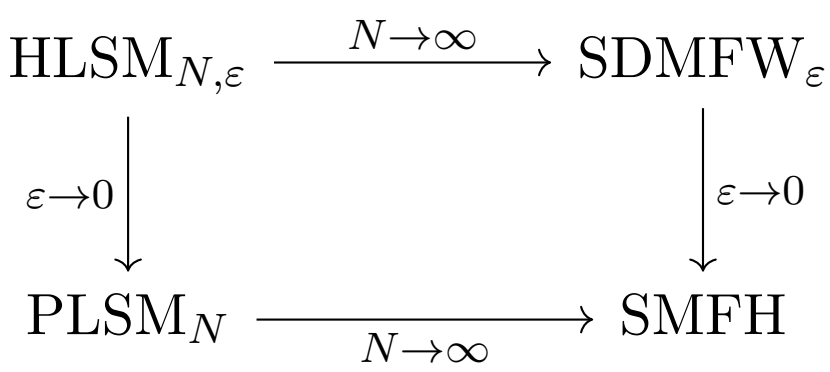
\begin{figure}[htbp]
    \centering
    \begin{tikzcd}[row sep=large, column sep=huge]
        \mathrm{HLSM}_{N,\varepsilon}
        \arrow[r, "N\to\infty"]
        \arrow[d, "\varepsilon\to0"']
        &
        \text{SDMFW}_{\varepsilon}
        \arrow[d, "\varepsilon\to0"]
        \\
        \mathrm{PLSM}_{N}
        \arrow[r, "N\to\infty"']
        &
        \text{SMFH}
    \end{tikzcd}
    \caption{Mean-field and Smoluchowski-Kramers limits of the
    stochastic quantization dynamics on $\mathbb{T}^2$.
    The arrows represent the respective limits under the
    hypotheses of \cite{czernik2026} and this work.}\label{commutative_diagram}
\end{figure}

Finally, we consider the Smoluchowski-Kramers approximation, which connects hyperbolic and parabolic stochastic quantization through the vanishing-inertia limit $\varepsilon\to0$. Building on \cite{Zine2025}'s work on singular stochastic wave equations on $\mathbb{T}^2$,  \cite{czernik2026} established this approximation for both the hyperbolic $O(N)$ linear sigma model and its mean-field limit. Together with the corresponding large-$N$ limits, their results yield the diagram in Figure~\ref{commutative_diagram}, with both orders of limits leading to the \eqref{eq:SMFH}.

The present work complements these dynamical results at the level of invariant measures. Under the stated renormalized second-moment condition, we identify the unique invariant law $\rho_\varepsilon$ of the \eqref{eq:SDMFW} and describe its $\varepsilon\to0$ behavior in Theorem~\ref{theorem:convergence_uniqueness_SMFH}.

\subsection{The SDMFW under an invariant measure}\label{sec:SDMFW_inv_msr}

We first fix the Wick-renormalization convention. Set
$\Lambda := 1-\Delta$ and let $\rho_0 := \mathcal{N}(0,\Lambda^{-1})$ be the massive ($m = 1$) Gaussian
free field measure. Let $\langle n\rangle=(1+|n|^2)^{1/2}$ denote the Japanese Bracket. We renormalize relative to $\rho_0$, defining
\[
    c_N := \mathbb{E}_{\rho_0}\bigl[(P_N u)(x)^2\bigr] = \sum_{\langle n\rangle\leq N}\langle n\rangle^{-2}
    \qquad \text{and} \qquad
    \wick{u^2} := \lim_{N\to\infty}\bigl((P_N u)^2-c_N\bigr),
\]
whenever the limit exists in a distributional sense that we specify later. This is the Wick ordering associated with the Euclidean $\Phi^4_2$-model as motivated previously. 

Suppose that the \eqref{eq:SDMFW} admits a stationary solution with law $\mu$, and that its renormalized second moment is well-defined. Stationarity implies, 
\[
    \mathbb{E}\bigl[\wick{u_\varepsilon(t)^2}\bigr]
    = \mathbb{E}\bigl[\wick{u_\varepsilon(0)^2}\bigr] = \mathbb{E}_{\mu}[\wick{u^2}] \eqcolon M_\mu.
\]
Thus, $M_\mu$ is a deterministic, time-independent spatial distribution. Once the stationary law is fixed, the mean-field term therefore becomes linear.

Writing $L_\mu := \Lambda+M_\mu$, we obtain
\begin{equation}\label{inv_SKA}
    \begin{cases}
        du_\eps=v_\eps\,dt,\\
        dv_\eps=-\eps^{-2}(L_\mu u_\eps+v_\eps)\,dt
        +\sqrt2\,\eps^{-2}\,dW_t,
    \end{cases}
\end{equation}
where $W$ is a cylindrical Wiener process on $L^2(\mathbb{T}^2)$.

The self-adjoint realization of $L_\mu$ and the weak interpretation of~\eqref{inv_SKA} are established in Section~\ref{Chpt_Invariant_Dynamics}.  
Equation~\eqref{inv_SKA} is linear once $M_\mu$ is fixed.  Projecting a
stationary solution onto any finite collection of eigenfunctions of
$L_\mu$ gives a finite-dimensional Ornstein-Uhlenbeck system.
The existence of an invariant probability measure for each such
projection forces every eigenvalue of $L_\mu$ to be strictly positive,
and the equation uniquely determines the invariant law. Proposition~\ref{prop_characterization_invmsr} yields
\[
    \mu = \mathcal{N}(0, L_\mu^{-1}) \otimes \mathcal{N}(0, \varepsilon^{-2}I) \in \mathcal P(\H^{-\sigma})
\]
Here the matrix denotes the distributional covariance.  For
$f\in H^\sigma(\T^2)$ and $g\in H^{1+\sigma}(\T^2)$,
\begin{align*}
    \EVMU{\bigl(\langle u,f\rangle+\langle v,g\rangle\bigr)^2}
    =
    \langle L_\mu^{-1}f,f\rangle_{L^2_x} +\varepsilon^{-2}\|g\|_{L^2_x}^2.
\end{align*}

\noindent
Hence, the natural support space is $\H$.
Indeed, the free stationary position and velocity fields belong almost surely to these respective spaces, while almost surely
neither belongs to the corresponding endpoint space $L^2(\mathbb{T}^2)$ or $H^{-1}(\mathbb{T}^2)$. This regularity is consistent with the constructions of singular stochastic wave dynamics in~\cite{GKOT2022,liu2025}.

\begin{remark}
    The present Wick ordering is compatible with a Da Prato-Debussche decomposition based Wick ordering used in \cite{GKOT2022,liu2025} who proved global pathwise existence of singular stochastic wave equations on $\T^2$ --- in the sense that their renormalization constant $c_{N,t}$ is of the same order as $c_N$.
    However, as $c_{N,t}$ is time dependent, subtracting this constant produces an explicitly time-dependent renormalization, for which stationarity of the solution law alone does not make the renormalized expectation constant in time, and therefore, in this instance, the renormalized equation would not behave as we expect. 
    Additionally, we note that \cite{Tolomeo2023Ergodicity} employed the same Wick ordering for his work on the conditional unique ergodictity of stochastic damped wave equations with polynomial nonlinearity on $\T^2$. 
\end{remark}

\subsection{Uniqueness through a self-consistency argument}
The law dependence of the \eqref{eq:SDMFW} prevents a direct application of standard uniqueness criteria for invariant measures as already noted in the study of the parabolic mean-field equation in \cite{Shen2022}.
Indeed, these criteria, such as the strong Feller Property combined with irreducibility (\cite{Da_Prato_Zabczyk_1996, HairerMattingly2018}) or Harris-type arguments (\cite{Hairer2011}), are typically formulated for a fixed Markov transition semigroup.
For a law-independent stochastic equation, with solution $X_t$, its transition probabilities
\[
    P_t(z,A):=\mathbb{P}(X_t^z\in A)
\]
are determined by the equation and do not depend on the initial distribution. In particular, the induced evolution of probability measures preserves convex combinations:
\[
    \bigl(a\nu_1+(1-a)\nu_2\bigr)P_t
    =
    a\nu_1P_t+(1-a)\nu_2P_t,
    \qquad 0\leq a\leq1.
\]
This remains true when the drift is nonlinear in the state.

The mean-field equation \eqref{eq:SDMFW} has a different
structure. Writing $\mathcal{S}_t\nu$ for the law at time $t$
of a solution initialized with law $\nu$, its drift involves
$
    M_{\mathcal{S}_t\nu}
    =
    \mathbb{E}_{\mathcal{S}_t\nu}[\wick{u^2}].
$
Thus, the coefficients themselves depend on the evolving law.
Under suitable well-posedness assumptions, the maps $\mathcal{S}_t$ form a semigroup on probability measures, but this semigroup is generally nonlinear and need not preserve convex combinations (see \cite{RehmeierRockner2025}). Consequently, it is not represented by a single law-independent Markov transition semigroup.

In the hyperbolic setting, the absence of spatial smoothing creates an additional obstacle to standard strong Feller arguments. \cite{Tolomeo2023Ergodicity} discusses this obstruction and the difficulties with the long-time estimates needed for asymptotic strong Feller methods, as proposed in \cite{Hairer2006}.
\cite{Tolomeo2023Ergodicity} overcomes these difficulties for the hyperbolic $P(\Phi)_2$ model through asymptotic strong Feller and coupling properties restricted to a group of regular shifts.
The argument yields ergodicity and conditional uniqueness for the fixed Markov semigroup associated with that equation.
For the \eqref{eq:SDMFW}, fixing a candidate invariant law $\mu$, freezes $M_\mu$ and similarly determines a time-homogeneous Markov semigroup. However, another candidate law $\nu$ may determine a different potential $M_\nu$ and therefore a different semigroup. Establishing uniqueness for each frozen equation does not compare these two laws. Rather, the question needed to be answered is which frozen potentials are consistent with their own invariant distributions.

For the parabolic mean-field equation \cite{Shen2022} prove uniqueness using decay estimates relative to the claimed unique invariant measure, the Gaussian free field, under a sufficiently large mass assumption. Our approach instead exploits the stationary covariance and does not require analogous long-time comparison estimates for the hyperbolic mean-field dynamics. 

A demonstrative simplification is the scalar self-consistency observation in Remark 5.5 in \cite{Shen2022}.
Suppose that $M_\mu$ is a spatial constant $c$. The Gaussian characterization, discussed in Section \ref{sec:SDMFW_inv_msr}, then gives the position covariance $(\Lambda+c)^{-1}$, where $\Lambda=1-\Delta$ and $c>-1$.
The self-consistent condition $c = \EV{u^2 - w^2}$ with $u  \sim \mathcal{N}(0,(\Lambda+c)^{-1})$ independent of $w \sim \mathcal{N}(0, \Lambda^{-1})$ then yields
\begin{align*}
    c
    =
    \sum_{n\in\mathbb{Z}^2}
    \left(
        \frac{1}{\langle n\rangle^2+c}
        -
        \frac{1}{\langle n\rangle^2}
    \right)
    =
    -c\sum_{n\in\mathbb{Z}^2}
    \frac{1}{
        \langle n\rangle^2
        (\langle n\rangle^2+c)
    }.
\end{align*}
The series converges and is strictly positive, so $c=0$. This identifies the only self-consistent constant potential.

For a general invariant law, however, $M_\mu$ need not be known a priori to be spatially constant.
Moreover, our assumptions initially allow $M_\mu$ to be a distribution. Thus, $L_\mu=\Lambda+M_\mu$ must first be given a rigorous operator realization.
Its inverse is generally not diagonal in the Fourier basis, and the scalar calculation above cannot be applied pointwise with $c$ replaced by $M_\mu(x)$.

Our argument extends self-consistency to this operator setting. We first construct $L_\mu$ through its quadratic form and show that the existence of an admissible invariant measure forces strict positivity. We then characterize that measure as Gaussian, with distributional covariance
\[
    \begin{pmatrix}
        L_\mu^{-1}&0\\
        0&\varepsilon^{-2}I
    \end{pmatrix}.
\]
The renormalized second moment is consequently determined by the covariance difference. Formally expanding $M_\mu$ in the Fourier basis $\{b_n\}_{n \in \Z^2}$, this gives
\[
    M_\mu
    =
    \sum_{n\in \Z^2} (L_\mu^{-1} - \Lambda^{-1}) (b_{n})(x) b_{-n}(x).
\]
The limit must be interpreted and justified carefully in an appropriate distribution space as the two unrenormalized covariance diagonals are not separately finite.

The resolvent identity
\[
    L_\mu^{-1}-\Lambda^{-1}
    =
    -L_\mu^{-1}M_\mu\Lambda^{-1}
\]
provides the operator analogue of the scalar identity. We use it to improve the regularity of $M_\mu$ and justify a positivity argument forcing $M_\mu=0$.
This resolves self-consistency for spatially varying potentials and establishes uniqueness in the class $\mathcal{W}$ without a large-mass assumption.

\subsection{Convergence of invariant measures}\label{sec:convergence_inv_msr}

We finally examine the stationary measure counterpart of the Smoluchowski-Kramers approximation in Figure~\ref{commutative_diagram}.
For the mean-field dynamics, the invariant laws identified above are
\[
    \rho_\varepsilon
    =
    \mathcal{N}(0,\Lambda^{-1})
    \otimes
    \mathcal{N}(0,\varepsilon^{-2}I),
    \qquad \text{and} \qquad
    \rho_0=\mathcal{N}(0,\Lambda^{-1}).
\]
Their uniqueness is understood within the classes specified in our theorems. 

The finite-$N$ hyperbolic Gibbs laws studied in \cite{liu2025} have an analogous Gaussian velocity factor, so the same arguments apply componentwise. We therefore present the argument for the mean-field law.

Fix $\sigma > 0$. Let $\pi_1$ and $\pi_2$ denote the coordinate projections
from $\mathcal{H}^{-\sigma}$ onto $H^{-\sigma}(\T^2)$ and $H^{-1-\sigma}(\T^2)$, respectively.
The position marginal satisfies
\[
    (\pi_1)_*\rho_\varepsilon=\rho_0
    \qquad\text{for every }\varepsilon>0.
\]
Thus, the stationary position distributions already agree with the parabolic invariant law. The singular behavior as $\varepsilon\to0$ occurs entirely in the velocity marginal
\[
    \eta_\varepsilon
    :=(\pi_2)_*\rho_\varepsilon
    =\mathcal{N}(0,\varepsilon^{-2}I).
\]

These velocity marginals do not converge weakly to a probability measure on $H^{-1-\sigma}(\T^2)$.
Indeed, choose a nonzero real-valued
$h\in C^\infty(\mathbb{T}^2)$ with $\|h\|_{L^2}=1$.
The continuous linear functional $v\mapsto\langle v,h\rangle$ has the one-dimensional distribution $\mathcal{N}(0,\varepsilon^{-2})$ under $\eta_\varepsilon$. Consequently, for every $R>0$,
\[
    \eta_\varepsilon
    \bigl(\{v:|\langle v,h\rangle|\leq R\}\bigr)
    =
    \mathbb{P}(|Z|\leq\varepsilon R)
    \longrightarrow0,
    \qquad Z\sim\mathcal{N}(0,1).
\]
Hence, the projected laws, and therefore $\eta_\varepsilon$ and $\rho_\varepsilon$, admit no weakly convergent subsequence as $\varepsilon\to0$ on their respective fixed state spaces.

To describe the velocity behavior, for a normed space $X$ we introduce the test-function class
\[
    C_{\mathrm{rv}}(X)
    :=
    \left\{
        f\in C_b(X):
        \lim_{R\to\infty}
        \sup_{\|x\|_X\geq R}|f(x)|=0
    \right\}.
\]
These functions vanish uniformly as the norm tends to infinity. The limiting statement considered here is
\begin{equation}\label{eq:stationary_SKA_marginals}
    (\pi_1)_*\rho_\varepsilon=\rho_0,
    \qquad
    \int_{H^{-1-\sigma}(\T^2)}f(v)\,\eta_\varepsilon(dv)
    \longrightarrow0
    \quad
    \text{for every }f\in C_{\mathrm{rv}}(H^{-1-\sigma}(\T^2)).
\end{equation}
The second assertion expresses escape of the velocity marginals from every bounded subset of $H^{-1-\sigma}(\T^2)$.
It does not assert weak convergence of $\rho_\eps$ to a probability measure.
Together with the exact identity for the position marginal, this describes the stationary measure Smoluchowski-Kramers limit.

\begin{remark}
    The more common notion of vague convergence, defined through continuous, compactly supported test functions, $C_c(H^{-1-\sigma}(\T^2))$, is trivial in this setting since $C_c(H^{-1-\sigma}(\T^2))=\{0\}$. A nonzero continuous function with compact support would force a closed ball of positive radius to be compact, which as a consequence of Riesz' Lemma is impossible in an infinite-dimensional normed space. Testing against $C_c(H^{-1-\sigma}(\T^2))$ would therefore give no information.
    By contrast, $C_{\mathrm{rv}}(H^{-1-\sigma}(\T^2))$ contains nonzero functions, such as $v\mapsto\exp(-\|v\|_{H^{-1-\sigma}(\T^2)}^2)$.
\end{remark}

\section{Uniqueness of the invariant dynamics of the SDMFW}\label{Chpt_Invariant_Dynamics}

Let $\sigma\in(0,\frac12)$ and let
\[
    \mu\in\mathcal P\bigl(\mathcal H^{-\sigma}(\T^2)\bigr)
    \qquad \text{with} \qquad
    \mathcal H^{-\sigma}
    =H^{-\sigma}(\T^2)\times H^{-1-\sigma}(\T^2),
\]
be an invariant law for \eqref{eq:SDMFW} such that
\[
    M_\mu:=\mathbb E_\mu[\wick{u^2}]
    \in H^{-\sigma}(\T^2).
\]
As explained in Section \ref{sec:SDMFW_inv_msr}, stationarity freezes
the mean-field term.  Hence, a stationary solution satisfies the linear
equation
\begin{equation}\label{eq:frozen_linear_dynamics}
    \begin{cases}
        du_\eps=v_\eps\,dt,\\
        dv_\eps=-\eps^{-2}(L_\mu u_\eps+v_\eps)\,dt
        +\sqrt2\,\eps^{-2}\,dW_t,
    \end{cases}
\end{equation}
with $\Lambda:=1-\Delta$ and $L_\mu:=\Lambda+M_\mu$. 

\noindent
We first construct $L_\mu$ as a self-adjoint operator with compact
resolvent that is bounded from below. We then prove sufficient Sobolev regularity and show that $L_\mu$ induces a Hilbert Space $\H_L$ with norm equivalent to the canonical norm of $\H^{-\sigma}$. Building on the framework established in \cite{DaPrato_Zabczyk_2014} to analyze stochastic differential equations (SDEs) in infinite dimensions, we prove that~\eqref{eq:frozen_linear_dynamics} admits a weak solution in Proposition~\ref{prop:mild} (in the sense of Theorem 5.4 in \cite{DaPrato_Zabczyk_2014}). Then, using the characterization Theorem 11.7 in~\cite{DaPrato_Zabczyk_2014} of invariant measures of linear infinite dimensional SDEs, we uniquely determine the invariant measure $\mu = \mathcal{N}(0, L_\mu^{-1}) \otimes \mathcal{N}(0, \varepsilon^{-2}I)$ of~\eqref{eq:frozen_linear_dynamics}. 
Finally, we are left to show that self-consistency forces $M_\mu = 0$ which implies that any invariant measure $\mu$ of the~\eqref{eq:SDMFW} equals $\rho_\eps = \mathcal{N}(0, \Lambda^{-1}) \otimes \mathcal{N}(0, \varepsilon^{-2}I)$.

\subsection{The operator \(L_{\mu}\) and its spectral basis}\label{sec_invdynamics_mildsolution}

For the remainder of this section we write $M = M_{\mu}$ and $L = L_{\mu}$.
The real-valued distribution $M$ defines on $H^1(\T^2)$ the symmetric
quadratic form
\[
    \mathfrak l[u,v]
    :=\langle \Lambda^{1/2}u,\Lambda^{1/2}v\rangle_{L^2_x}
      +\langle M,u\overline v\rangle,
\]
where the second term is understood by Sobolev duality.

This first lemma utilizes the KLMN Theorem to prove that the operator $L$ can be realized as a self-adjoint operator on $L^2(\T^2)$.
\begin{lemma}\label{lem:relative_bound}
Suppose $M\in H^{-\sigma}(\T^2)$ for some $\sigma\in(0,1)$.
For every $\delta>0$ there exists $C_\delta>0$ such that
\[
    |\langle M,|u|^2\rangle|
    \leq
    \delta\langle\Lambda u,u\rangle_{L^2_x}
    +C_\delta\|u\|_{L^2_x}^2,
    \qquad u\in H^1(\T^2).
\]
Consequently, $\mathfrak l$ is closed and bounded from below on
$H^1(\T^2)$ and determines a self-adjoint operator $L$ on
$L^2(\T^2)$ with form domain $H^1(\T^2)$.
\end{lemma}

\begin{proof}
Choose $\eta>0$ such that $2\eta+\sigma\leq1$.  By Sobolev duality,
the negative-order product estimate in \cite[Proposition~2]{Benyi2025},
and Sobolev embedding,
\begin{align*}
    |\langle M,|u|^2\rangle|
    \leq
    \|Mu\|_{H^{-1+\eta}_x}\|u\|_{H^{1-\eta}_x}
    \lesssim
    \|M\|_{H^{-\sigma}_x}\|u\|_{H^{1-\eta}_x}^2.
\end{align*}
Interpolating between fractional Sobolev spaces, for any $\tilde{\delta} > 0$, there exists a constant $C(\tilde{\delta}) > 0$ such that
\[
    \|u\|_{H^{1-\eta}_x}^2
    \leq
    \widetilde\delta\|u\|_{H^1_x}^2
    +C(\widetilde{\delta})\|u\|_{L^2_x}^2.
\]
Taking $\widetilde\delta$ sufficiently small proves the stated bound.
The final assertion follows from the KLMN theorem
\cite[Theorem~10.21]{schmudgen_operators}.
\end{proof}

The next lemma demonstrates that $L$ has compact resolvent and characterizes its spectrum in relation to the operator $\Lambda$. 
\begin{lemma}\label{lemma_equ_ev}
Suppose $M\in H^{-\sigma}(\T^2)$ for some $\sigma\in(0,1)$. There exists $c>0$ such that $L_c:=L+c$ is strictly positive and
\begin{equation}\label{eq:form_equivalence}
    C_1\langle\Lambda u,u\rangle_{L^2_x}
    \leq
    \langle L_cu,u\rangle_{L^2_x}
    \leq
    C_2\langle\Lambda u,u\rangle_{L^2_x},
    \qquad u\in H^1(\T^2),
\end{equation}
for some $C_1,C_2>0$.  Moreover, $L$ has compact resolvent.  In
particular, there exist a real orthonormal basis $\{e_k\}_{k\in \N}$ of
$L^2(\T^2)$ and real eigenvalues $(\lambda_k)_{k\geq1}$, repeated
according to multiplicity and satisfying $\lambda_k\to+\infty$, such
that
\[
    Le_k=\lambda_ke_k.
\]
If $(\nu_k)_{k\geq1}$ are the eigenvalues of $L_c$ and
$(\ell_k)_{k\geq1}$ those of $\Lambda$, both in nondecreasing order,
then
\[
    C_1\ell_k\leq\nu_k\leq C_2\ell_k,
    \qquad k\geq1.
\]
\end{lemma}

\begin{proof}
By Lemma \ref{lem:relative_bound}, $L$ is bounded from below.  Choose
$c>0$ so that $L_c$ is strictly positive.  The relative form bound,
combined with strict positivity of $L_c$, gives
\eqref{eq:form_equivalence}.

For $f\in L^2(\T^2)$ and $u=L_c^{-1}f$, the lower bound in
\eqref{eq:form_equivalence} gives
\[
    \|u\|_{H^1_x}^2
    \lesssim
    \langle L_cu,u\rangle_{L^2_x}
    =\langle f,u\rangle_{L^2_x}
    \lesssim
    \|f\|_{L^2_x}\|u\|_{H^1_x}.
\]
Thus, $L_c^{-1}:L^2(\T^2)\to H^1(\T^2)$ is bounded.  Since
$H^1(\T^2)\hookrightarrow L^2(\T^2)$ is compact, $L_c^{-1}$ is
compact on $L^2$; hence both $L_c$ and $L$ have compact resolvent.
The spectral assertions follow from the spectral theorem.  Finally,
the eigenvalue comparison follows from \eqref{eq:form_equivalence} and
the min-max principle.
\end{proof}
Since the previous Lemma establishes that for some $c > 0$, $L_c \coloneq L + c$ is a strictly positive self-adjoint operator, we can define $\H_L$ as the closure of $L^2(\T^2) \times L^2(\T^2)$ with respect to the norm induced by the inner product
\begin{align}\label{eq:inner_product}
    \langle \mathbf{u}, \mathbf{v} \rangle_{L_c} \coloneqq \langle L_c^{-\sigma} u_1, v_1 \rangle_{L^2_x} + \langle L_c^{-1-\sigma} u_2, v_2 \rangle_{L^2_x},
\end{align}
with respect to the normalized inner product on $L^2(\T^2)$:
\begin{align*}
    \langle u, v \rangle_{L^2_x} \coloneq (2\pi)^{-2} \int_{\T^2} u(x) \overline{v(x)} \, dx.
\end{align*}

Finally, the following lemma is the point at which the restriction $\sigma<\frac12$ enters the spectral projection argument. The Corollary following this Lemma, infers that the Hilbert spaces $\H_L$ and $\H^{-\sigma}$ are equivalent. 
\begin{lemma}\label{lem:eigenfunction_regularity}
Let $\sigma\in(0,\frac12)$ and $M\in H^{-\sigma}(\T^2)$.  For the
constant $c$ in Lemma \ref{lemma_equ_ev}, the operator $L_c$ extends
to an isomorphism
\begin{equation}\label{eq:Lc_sobolev_isomorphism}
    L_c:H^{1+\sigma}(\T^2)
    \longrightarrow H^{-1+\sigma}(\T^2).
\end{equation}
Consequently,
\[
    e_k\in H^{1+\sigma}(\T^2),
    \qquad k\geq1,
\]
and the algebraic span
\[
    \mathscr E_L:=\operatorname{span}\{e_k:k\geq1\}
\]
is dense in both $H^{1+\sigma}(\T^2)$ and
$H^\sigma(\T^2)$.
\end{lemma}

\begin{proof}
From the product estimate \cite{Benyi2025} Proposition 2 and Sobolev's inequality if follows that multiplication by an element of
$H^{1+\sigma}(\T^2)$ is bounded on $H^{-\sigma}(\T^2)$.  Hence,
\begin{equation}\label{eq:potential_multiplier_bound}
    \|Mf\|_{H^{-\sigma}_x}
    \lesssim
    \|M\|_{H^{-\sigma}_x}\|f\|_{H^{1+\sigma}_x}.
\end{equation}
Because $0 < \sigma<\frac12$, the embedding
\[
    H^{-\sigma}(\T^2)
    \hookrightarrow H^{-1+\sigma}(\T^2)
\]
is compact.  It follows that multiplication by $M$ is a compact
operator from $H^{1+\sigma}(\T^2)$ to $H^{-1+\sigma}(\T^2)$.  The same is true of
the zeroth-order term $cI$.  Since
\[
    \Lambda:H^{1+\sigma}(\T^2)
    \longrightarrow H^{-1+\sigma}(\T^2)
\]
is an isomorphism, $K = \Lambda^{-1}(M+cI)$ is compact. Additionally, $L_c$ has trivial kernel: if $L_cf=0$, then
$f\in H^{1+\sigma}(\T^2)\subset H^1(\T^2)$, and strict positivity of the form of $L_c$ gives $f=0$.
Thus, writing $L_c=\Lambda+M+c = \Lambda(I + K)$, the Fredholm alternative \cite[Theorem~6.6]{Brezis2011} implies that $L_c$ is Fredholm of index zero between
these spaces. Finally, from the bounded inverse theorem follows \eqref{eq:Lc_sobolev_isomorphism}.

Since $(\lambda_k+c)e_k\in L^2(\T^2)\subset H^{-1+\sigma}(\T^2)$,
\eqref{eq:Lc_sobolev_isomorphism} gives $e_k\in H^{1+\sigma}(\T^2)$.
It remains to prove density.  Let $f\in H^{1+\sigma}(\T^2)$ and set
$g=L_cf\in H^{-1+\sigma}(\T^2)$.  The finite linear span of $\{e_k\}_{k \in \N}$ is
dense in $H^{-1+\sigma}(\T^2)$: indeed, $L^2(\T^2)$ is dense in
$H^{-1+\sigma}(\T^2)$, while the finite eigenfunction expansions are dense
in $L^2(\T^2)$.  Choose $g_N\in\mathscr E_L$ with
$g_N\to g$ in $H^{-1+\sigma}(\T^2)$ and set $f_N=L_c^{-1}g_N$.  Then
$f_N\in\mathscr E_L$ and, by \eqref{eq:Lc_sobolev_isomorphism},
$f_N\to f$ in $H^{1+\sigma}(\T^2)$.  Thus, $\mathscr E_L$ is dense in
$H^{1+\sigma}(\T^2)$.  Its density in $H^\sigma(\T^2)$ follows from the density
of $H^{1+\sigma}(\T^2)$ in $H^\sigma(\T^2)$.
\end{proof}
For \(s\geq0\), equip
\[
 H_{L_c}^s:=D(L_c^{s/2})
\]
with the norm 
\[
\|f\|_{H_{L_c}^s}:=\|L_c^{s/2}f\|_{L^2_x}.
\]
For \(s<0\), let \(H_{L_c}^s\) be the completion of \(L^2(\T^2)\)
in this norm, equivalently the \(L^2\)-pivot dual of \(H_{L_c}^{-s}\).

\begin{corollary}[Identification of the Sobolev phase space]\label{cor:equ_norms}
For every \(s\in[-1-\sigma,1+\sigma]\),
\begin{align*}
    H_{L_c}^s=H^s(\T^2)
\end{align*}
with equivalent norms. In particular,
\begin{align*}
 \H_L=H_{L_c}^{-\sigma}\times H_{L_c}^{-1-\sigma}
     =\H^{-\sigma}.
\end{align*}
These identifications preserve the Borel sigma-algebras.
\end{corollary}

\begin{proof}

Form equivalence of $L_c$ and $\Lambda$, established in Lemma~\ref{lemma_equ_ev}, identifies \(H_{L_c}^1\) with \(H^1(\T^2)\). By the Heinz inequality for self-adjoint operators (\cite[Proposition 10.14]{schmudgen_operators}) form equivalence of $L_c$ and $\Lambda$ translates to the form equivalence of $L_c^{s}$ and $\Lambda^{s}$ for $s \in [0,1]$. Then, Corollary 10.12 in \cite{schmudgen_operators} implies that $L_c^{-s}$ and $\Lambda^{-s}$ are form equivalent. In particular, for \(g\in L^2(\T^2)\),
\[
 \|g\|_{H^{-1+\sigma}}
 \asymp\|L_c^{-(1-\sigma)/2}g\|_{L^2}.
\]
For \(f\in\mathscr E_L\), \eqref{eq:Lc_sobolev_isomorphism} now gives
\[
 \|f\|_{H^{1+\sigma}}
 \asymp\|L_c f\|_{H^{-1+\sigma}}
 \asymp\|L_c^{(1+\sigma)/2}f\|_{L^2}.
\]
By Lemma~\ref{lem:eigenfunction_regularity}, \(\mathscr E_L\) is dense
in \(H^{1+\sigma}(\T^2)\), and it is dense in \(H_{L_c}^{1+\sigma}\) by spectral
truncation.  Both spaces embed continuously into \(L^2(\T^2)\), so completion identifies
them as the same space. Interpolation and duality give the remaining indices.
The assertion about Borel sets follows from equivalence of norms.
\end{proof}

\subsection{Weak solution of the frozen SDE and its invariant measure}\label{sec_invdynamics_unique}
Throughout this subsection, all Hilbert spaces are real. We work on
\(\H:=\H_L\), equipped with the inner product
\eqref{eq:inner_product}, and set \(U:=L^2(\T^2)\).
By Corollary~\ref{cor:equ_norms}, \(\H\) is the original Sobolev space \(\H^{-\sigma}\) with an equivalent norm.
Define the operators 
\begin{equation}\label{eq:generator}
 \begin{split}
 \mathcal A(u,v)&=(v,-\eps^{-2}(Lu+v)),\\
 B_\eps h&=(0,\sqrt2\,\eps^{-2}h)
 \end{split}
\end{equation}
with domains $D(\mathcal A)=H_{L_c}^{1-\sigma}\times H_{L_c}^{-\sigma} =H^{1-\sigma}(\T^2)\times H^{-\sigma}(\T^2)$ and $D(B) = U$.
Here \(L=L_c-cI\) is extended continuously from
\(H_{L_c}^{1-\sigma}\) to \(H_{L_c}^{-1-\sigma}\).
This extension agrees with the form realization on their common
domain and does not require a definition of \(Mu\) for arbitrary
\(u\in H^{-\sigma}\).

We establish that~\eqref{eq:frozen_linear_dynamics} admits a weak solution in the sense of \cite{DaPrato_Zabczyk_2014} Sections 5.1--5.2. More precisely, on a filtered probability space, with a cylindrical \(U\)-Wiener process \(W\) relative
to the filtration, an adapted progressively measurable
\(\H\)-valued process \(X\) is a weak solution if its paths are locally
integrable in \(\H\) and, for every \(h\in D(\mathcal A^*)\),
\begin{equation}\label{eq:weak_frozen}
 \langle X(t),h\rangle_{\mathcal{H}}
 =\langle X(0),h\rangle_{\mathcal H}
  +\int_0^t\langle X(s),\mathcal A^*h\rangle_{\mathcal{H}}\,ds
  +\int_0^t\langle B_\eps^*h,dW_s\rangle_U .
\end{equation}
The identity is understood almost surely for each \(t\). We use
the continuous version obtained below. All adjoints in this
subsection are taken with respect to the fixed inner product of
\(\H\), unless explicitly stated otherwise.

\begin{proposition}[Existence and uniqueness of a weak solution]\label{prop:mild}
Suppose $M\in H^{-\sigma}(\T^2)$ for some $\sigma\in(0,\frac12)$.
The operator \(\mathcal A\) generates a strongly continuous semigroup \(S(t)\) for \(t\geq0\) on \(\H\).
For every \(\H\)-valued \(\mathcal F_0\)-measurable random variable
\(X_0\), the frozen equation has a unique weak solution, up to
modification, given by
\begin{equation}\label{eq:mild}
 X(t)=S(t)X_0+\int_0^t S(t-s)B_\eps\,dW_s
\end{equation}
for $t \geq 0$.

\noindent
\end{proposition}

\begin{proof} 
We show that the stochastic process $X(t)$ satisfies the conditions of Theorem 5.4 in \cite{DaPrato_Zabczyk_2014}. Then existence and uniqueness of a weak solution follow from this theorem. We are left to establish that $\mathcal A$ generates a strongly continuous semigroup in $\H$ and that the stochastic convolution is well-defined.
Decompose \(\mathcal A=\mathcal J+\mathcal R\), where
\[
 \mathcal J(u,v)=(v,-\eps^{-2}L_c u),\qquad
 \mathcal R(u,v)=(0,\eps^{-2}cu-\eps^{-2}v).
\]
On finite spectral sums, define, for every \(t\in\R\),
\begin{equation}\label{eq:abstract_wave_group}
 G(t)=
 \begin{pmatrix}
 \cos(tL_c^{1/2}/\eps)&
       \eps L_c^{-1/2}\sin(tL_c^{1/2}/\eps)\\
 -\eps^{-1}L_c^{1/2}\sin(tL_c^{1/2}/\eps)&
       \cos(tL_c^{1/2}/\eps)
 \end{pmatrix}.
\end{equation}
The spectral calculus interprets the off-diagonal entries between
the two components of \(\H\). On the \(k\)-th mode this is the
solution operator for
\(u_k'=v_k,\ v_k'=-\eps^{-2}\nu_k u_k\), where
\(\nu_k=\lambda_k+c\) are the eigenvalues of $L_c$. It preserves the equivalent norm
\begin{equation}\label{eq:equ_norm}
 \|(u,v)\|_{\H,\eps}^2
 =\sum_{k\geq1}\left(\nu_k^{-\sigma}|u_k|^2
             +\eps^2\nu_k^{-1-\sigma}|v_k|^2\right).
\end{equation}
Consequently, \(G(t)\) extends to an isometry for the norm
\(\|\cdot\|_{\H,\eps}\).
The trigonometric addition formulas give \(G(t+s)=G(t)G(s)\)
and \(G(-t)=G(t)^{-1}\). Strong continuity follows by approximation
with finite spectral sums.

Its generator is \(\mathcal J\). Indeed, differentiation in each
mode gives the stated action, and the condition
\(\mathcal J(u,v)\in \H\) is exactly
\[
 \sum_k\nu_k^{1-\sigma}|u_k|^2
       +\sum_k\nu_k^{-\sigma}|v_k|^2<\infty.
\]
For elements satisfying this condition, the spectral difference
quotients converge in \(\H\), by dominated convergence; conversely,
convergence in \(\H\) forces this condition by coordinatewise
differentiation. Thus, \(D(\mathcal J)=D(\mathcal A)\).

The embedding \(H_{L_c}^{-\sigma}\hookrightarrow H_{L_c}^{-1-\sigma}\)
makes \(\mathcal R\) bounded on \(\H\). Applying the bounded
perturbation theorem to \(\mathcal J+\mathcal R\) shows that \(\mathcal A\) generates a strongly continuous semigroups.

To define the stochastic convolution as an element of
$L^2(\Omega;\mathcal H)$, it suffices to verify that
\[
    \int_0^t
    \|S(t-s)B_\eps\|_{\mathscr L_2(U,\mathcal H)}^2\,ds
    <\infty
    \qquad\text{for every }t>0.
\]
Here $\mathscr L_2(U,\mathcal H)$ denotes the space of
Hilbert--Schmidt operators from $U$ to $\mathcal H$, with
\[
    \|T\|_{\mathscr L_2(U,\mathcal H)}^2
    :=\sum_{k\geq1}\|Te_k\|_{\mathcal H}^2,
\]
where $\{e_k\}_{k \in \N}$ is any orthonormal basis of $U$.
The Hilbert space Itô isometry then gives
\[
    \mathbb E\left\|
        \int_0^t S(t-s)B_\eps\,dW_s
    \right\|_{\mathcal H}^2
    =
    \int_0^t
    \|S(t-s)B_\eps\|_{\mathscr L_2(U,\mathcal H)}^2\,ds.
\]
For that we return to the original inner
product \eqref{eq:inner_product}, without the additional factor
\(\eps^2\) in \eqref{eq:equ_norm}. In this norm,
\[
 \|B_\eps\|_{\mathscr L_2(U,\H)}^2
 =2\eps^{-4}\sum_{k\geq1}\nu_k^{-1-\sigma}<\infty,
\]
by Lemma~\ref{lemma_equ_ev} and \(\sum_{n\in\Z^2}(1+|n|^2)^{-1-\sigma}<\infty\).
Therefore, for every \(t>0\),
\[
 \int_0^t\|S(t-s)B_\eps\|_{\mathscr L_2(U,\H)}^2\,ds
 \leq t\sup_{0\leq r\leq t}\|S(r)\|_{\mathscr L(\H)}^2
             \|B_\eps\|_{\mathscr L_2(U,\H)}^2<\infty.
\]
This finishes the proof.
\end{proof}

\begin{remark}
The shift \(c\) is used to construct equivalent norms, not to change
the drift. The operator in \eqref{eq:generator} contains \(L\),
not \(L+c\). Different admissible shifts yield the same realization:
they have equivalent Hilbert norms and the same action on the
finite spectral core. In particular, this construction does not
presuppose the positivity of \(L\) that is later deduced from
stationarity.
\end{remark}

\begin{proposition}[Characterization of the invariant measure]\label{prop_characterization_invmsr}
Fix \(0<\sigma<\frac12\), \(\eps>0\), and a real
\(M\in H^{-\sigma}(\T^2)\), and let \(L=\Lambda+M\) be its form
realization. Suppose \eqref{eq:frozen_linear_dynamics} has an invariant probability measure
\(\eta\in\mathcal P(\H^{-\sigma})\). Then:
\begin{enumerate}
 \item \emph{Positivity:} \(L\geq mI\) for some \(m>0\), and the forms of \(L\) and
       \(\Lambda\) are equivalent on \(H^1(\T^2)\).
 \item \emph{Exponential decay:} There exist \(C\geq1\) and \(\omega>0\) such that the operator norm in $\H$ of the semigroup $S(t)$ has exponential decay:
 \begin{equation}\label{etA_exponential_decay}
    \|S(t)\|_{\H}\leq Ce^{-\omega t},
       \qquad t\geq0.
 \end{equation}
 \item \emph{Gaussianity and uniqueness.} The invariant measure is unique and centered Gaussian.
       Its distributional covariance is
 \begin{equation}\label{eq:characterized_covariance}
       \begin{pmatrix}L^{-1}&0\\0&\eps^{-2}I\end{pmatrix}.
 \end{equation}
 In particular, for real \(f\in H^\sigma(\T^2)\) and
 \(g\in H^{1+\sigma}(\T^2)\),
 \begin{equation}\label{eq:characteristic_functional_mu}
  \int_{\mathcal H} e^{i\langle u,f\rangle+i\langle v,g\rangle}\,d\eta(u,v)
   =\exp\!\left[-\frac12\left(
       \langle L^{-1}f,f\rangle_{L^2_x}
          +\eps^{-2}\|g\|_{L^2_x}^2\right)\right]
 \end{equation}
 where the term on the left is understood by Sobolev duality.
\end{enumerate}
\end{proposition}

\begin{proof}
Set
\[
 Q_t=\int_0^t S(r)B_\eps B_\eps^*S(r)^*\,dr.
\]

\emph{Positivity of \(L\):}
For the positivity argument, choose $c\geq0$ such that
$L_c=L+cI$ is strictly positive, and equip $\H=\H_L$
with the inner product \eqref{eq:inner_product}.
Adjoints and traces in this step refer to this inner product.
The characterization theorem for invariant measures of linear
stochastic equations Theorem 11.7 in \cite{DaPrato_Zabczyk_2014}
implies, from existence of an invariant probability measure, that
\begin{equation}\label{eq:uniform_trace_bound}
 \sup_{t\geq0}\operatorname{Tr}_{\mathcal H} Q_t<\infty.
\end{equation}
Let $\{e_k\}_{k \in \N}$ be the spectral basis of $L$ with corresponding eigenvalues $\{\lambda_k\}_{k \in \N}$ then by the Hilbert--Schmidt trace identity,
\[
 \operatorname{Tr}_{\mathcal H} Q_t
   =\sum_{k\in\N}\int_0^t\|S(r)B_\eps e_k\|_{\mathcal H}^2\,dr.
\]
Suppose \(\lambda_k\leq0\). Write
\(S(t)B_\eps e_k=(a_k(t)e_k,a_k'(t)e_k)\). Then
\[
 \eps^2a_k''+a_k'+\lambda_ka_k=0,\qquad
 a_k(0)=0,\quad a_k'(0)=\sqrt2\,\eps^{-2}.
\]
With
\[
 r_\pm=\frac{-1\pm\sqrt{1-4\eps^2\lambda_k}}{2\eps^2},
\]
these distinct real roots give
\[
 a_k(t)=\frac{\sqrt2}{\eps^2(r_+-r_-)}
              \left(e^{r_+t}-e^{r_-t}\right).
\]
If \(\lambda_k=0\), this equals
\(\sqrt2(1-e^{-t/\eps^2})\); if \(\lambda_k<0\), it grows
exponentially. In either case,
\[
 \operatorname{Tr}_{\mathcal{H}} Q_t
 \geq (\lambda_k + c)^{-\sigma}\int_0^t|a_k(r)|^2\,dr
 \longrightarrow\infty,
\]
contrary to \eqref{eq:uniform_trace_bound}.
Thus, every \(\lambda_k>0\). Since \(L\) is bounded below and has
compact resolvent, \(m:=\inf_k\lambda_k>0\).
Moreover,
\[
 \frac{m}{m+c}L_c\leq L\leq L_c
\]
in the form sense. Lemma~\ref{lemma_equ_ev} proves the asserted
form equivalence.

\emph{Exponential decay:}
Since $L\geq mI$, we can choose $c=0$. 
To prove exponential decay of the semigroup we consider the $\eps$-adjusted energy functional
\[
 \mathcal E(u,v)
   :=\|L^{-\sigma/2}u\|_{L^2_x}^2
        +\eps^2\|L^{-(1+\sigma)/2}v\|_{L^2_x}^2
\]
that is a squared norm equivalent to \eqref{eq:inner_product} on $\H$. 

We first establish the energy estimate on finite spectral
subspaces. For $N\geq1$, set
\[
 \H_N:=\operatorname{span}\{(e_k,0),(0,e_k):1\leq k\leq N\}.
\]
The block-diagonal spectral representation of the generator (see the proof of Proposition~\ref{prop:mild})
implies that $\H_N$ is invariant under $S(t)$. Hence, for
$X_0\in \H_N$, the deterministic solution
\[
 X(t)=S(t)X_0=(u(t),v(t))
\]
has the finite expansions
\[
 u(t)=\sum_{k=1}^N u_k(t)e_k,
 \qquad
 v(t)=\sum_{k=1}^N v_k(t)e_k,
\]
where
\[
 u_k'(t)=v_k(t),\qquad
 \eps^2v_k'(t)=-\lambda_k u_k(t)-v_k(t).
\]
We refer to such a solution as a finite spectral solution.
All subsequent differentiations are therefore ordinary
differentiations of finite sums. Crucially, the estimates
below are independent of $N$ and of $X_0$.

Put
\[
 w(t)=L^{-(1+\sigma)/2}u(t),\qquad
 z(t)=L^{-(1+\sigma)/2}v(t).
\]
Then
\[
 w'(t)=z(t),\qquad \eps^2z'(t)=-Lw(t)-z(t),
\]
and the energy of the finite spectral solution is
\[
 \mathcal E(X(t)) 
 =\|L^{1/2}w(t)\|_{L^2_x}^2
    +\eps^2\|z(t)\|_{L^2_x}^2.
\]
For \(\delta>0\), define
\[
 \mathcal V(t)
   :=\mathcal E (X(t)) +2\delta\eps^2\langle w(t),z(t)\rangle_{L^2_x}
                    +\delta\|w(t)\|_{L^2_x}^2.
\]
Choose \(\delta\eps^2\leq\frac12\) and
\(\delta\eps/\sqrt m\leq\frac12\).
Cauchy--Schwarz and \(L\geq mI\) give
\[
 \frac12\mathcal E(X(t))\leq\mathcal V(t)
       \leq\left(\frac32+\frac{\delta}{m}\right)\mathcal E(X(t)).
\]
Direct differentiation yields
\[
 \mathcal V'(t)
  =-2(1-\delta\eps^2)\|z(t)\|_{L^2_x}^2
      -2\delta\|L^{1/2}w(t)\|_{L^2_x}^2
  \leq-\min\{\eps^{-2},2\delta\}\mathcal E(X(t)).
\]

Set
\[
 K:=\frac32+\frac{\delta}{m},
 \qquad q:=\min\{\eps^{-2},2\delta\}>0.
\]
The preceding estimates imply
\[
 \mathcal V'(t)
 \leq-q\mathcal E(X(t))
 \leq-\frac{q}{K}\mathcal V(t).
\]
Gronwall's inequality therefore yields
\[
 \mathcal E(S(t)X_0)
 \leq 2\mathcal V(t)
 \leq 2K e^{-qt/K}\mathcal E(X_0).
\]
By equivalence of norms,
\[
 \|S(t)X_0\|_{\mathcal H}
 \les \sqrt{2K} e^{-qt/(2K)}\|X_0\|_{\mathcal H}
\]
for $X_0\in\bigcup_{N\geq1}\H_N$.

\noindent
The constants are independent of the spectral cutoff.
Since $\bigcup_{N\geq1}\H_N$ is dense in $\H$ and $S(t)$
is bounded for each fixed $t\geq0$, this estimate extends
to every $X_0\in\H$. Taking the supremum over
$\|X_0\|_{\mathcal H}=1$ proves
\[
 \|S(t)\|_{\mathcal H}
 \leq C e^{-\omega t},
\]
for some $C\geq1$ and $\omega>0$.

\emph{Gaussianity and uniqueness:}
Utilizing the exponential decay we can apply Proposition 11.10 in \cite{DaPrato_Zabczyk_2014} to obtain that $\eta$ is the unique invariant measure of \eqref{eq:frozen_linear_dynamics}. Then, Theorem 11.7 in \cite{DaPrato_Zabczyk_2014} implies that $\eta = \mathcal{N}(0, Q_\infty)$ with $Q_\infty = \lim_{t\to\infty} Q_t$.
Finally, it remains to explicitly characterize $Q_\infty$. 
The vectors
\[
 F_k=(\lambda_k^{\sigma/2}e_k,0),\qquad
 G_k=(0,\lambda_k^{(1+\sigma)/2}e_k)
\]
form an orthonormal basis of \(\H\). The orthogonal planes
\(\operatorname{span}\{F_k,G_k\}\) reduce both \(S(t)\) and
\(B_\eps B_\eps^*\). On the \(k\)-th plane, the drift and diffusion
matrices in this orthonormal basis are
\[
 \widetilde A_k=
 \begin{pmatrix}
 0&\lambda_k^{1/2}\\
 -\eps^{-2}\lambda_k^{1/2}&-\eps^{-2}
 \end{pmatrix},\qquad
 \widetilde B_k=
 \begin{pmatrix}0\\ \sqrt2\,\eps^{-2}\lambda_k^{-(1+\sigma)/2}\end{pmatrix}.
\]
Thus, the covariance block \(C_k\) solves the usual Euclidean
Lyapunov equation
\begin{equation}\label{lyapunov_eq}
 \widetilde A_kC_k+C_k\widetilde A_k^{\mathsf T}
       +\widetilde B_k\widetilde B_k^{\mathsf T}=0.
\end{equation}
A direct substitution gives
\[
 C_k=\begin{pmatrix}
       \lambda_k^{-1-\sigma}&0\\
       0&\eps^{-2}\lambda_k^{-1-\sigma}
     \end{pmatrix}.
\]
There are no cross-mode covariance blocks, because the covariance
integrand is block diagonal. Consequently, the covariance operator
relative to the inner product of \(\H\) is
\begin{equation}\label{eq:hilbert_covariance}
 Q_\infty=
 \begin{pmatrix}
    L^{-1-\sigma}&0\\
    0&\eps^{-2}L^{-1-\sigma}
 \end{pmatrix}.
\end{equation}
For \(X=(u,v)\sim\mathcal N(0,Q_\infty)\), returning from the
orthonormal coordinates to the distributional coordinates yields
\[
 \mathbb E[\langle u,e_k\rangle\langle u,e_j\rangle]
     =\delta_{kj}\lambda_k^{-1},\qquad
 \mathbb E[\langle v,e_k\rangle\langle v,e_j\rangle]
     =\delta_{kj}\eps^{-2},
 \qquad
 \mathbb E[\langle u,e_k\rangle\langle v,e_j\rangle]=0.
\]
This proves \eqref{eq:characteristic_functional_mu} for finite
spectral test functions. Their density in
\(H^\sigma(\T^2)\times H^{1+\sigma}(\T^2)\), established in Lemma~\ref{lem:eigenfunction_regularity}, together with bounded convergence
on the left and continuity of the quadratic form on the right,
proves the identity for all stated test functions.
It is \eqref{eq:characterized_covariance}, rather than
\eqref{eq:hilbert_covariance}, that denotes the distributional
covariance of the invariant field.
\end{proof}

\subsection{Uniqueness by self-consistency}
Proposition \ref{prop_characterization_invmsr} reduces the remaining
argument to a self-consistency problem.  Namely, every admissible invariant law
has position marginal $\mathcal N(0,L^{-1})$. Utilizing this self-consistency relation we show that $M_\mu = 0$. This directly implies that $\Lambda = L$ which proves uniqueness. 
The following lemma identifies this
renormalized covariance diagonal and improves the regularity of
$M_\mu$ to $L^2(\T^2)$.

\begin{lemma}\label{lem:renormalized_variance_convergence}
    Let $\sigma \in (0, \frac 12)$ and $u \sim \mathcal{N}(0, L^{-1})$ be a centered Gaussian random distribution. Assume that $M_\mu \in H^{-\sigma}(\T^2)$ and set $L = \Lambda + M_\mu$.

    Set
    \[
        \Gamma_N
        \coloneqq
        \{n\in\Z^2:\langle n\rangle\leq N\},
        \qquad
        \P_Nf
        \coloneqq
        \sum_{n\in\Gamma_N}
        \langle f,b_n\rangle_{L^2_x}b_n,
    \]
    where \(b_n(x)=e^{in\cdot x}\), and define
    \[
        M_N
        \coloneqq
        \mathbb E_\mu[(\P_Nu)^2]
        -
        \sum_{n\in\Gamma_N}\langle n\rangle^{-2}.
    \]
    Then $M_\mu\in L^2(\T^2)$ and
    $M_N\to M_\mu$ strongly in $L^2(\T^2)$.
    More precisely, writing $D=L^{-1}-\Lambda^{-1}$, we have
    \begin{equation}\label{eq:M_mu_expansion}
        M_\mu
        =
        \lim_{N\to\infty}
        \sum_{n\in\Gamma_N}\P_N(Db_n)b_{-n}
        =
        \lim_{N\to\infty}
        \sum_{n\in\Gamma_N}(Db_n)b_{-n},
    \end{equation}
    where both limits are taken in $L^2(\T^2)$.
\end{lemma}

\begin{proof}
    Since the covariance form of \(u\) is \(L^{-1}\), we have
    \[
        \mathbb E_\mu
        \left[
            \langle u,b_n\rangle_{L^2_x}
            \langle u,b_m\rangle_{L^2_x}
        \right]
        =
        \langle L^{-1}b_{-n},b_m\rangle_{L^2_x}.
    \]
    Consequently, after relabeling \(n\mapsto-n\),
    \[
        \mathbb E_\mu[(\P_Nu)^2]
        =
        \sum_{n\in\Gamma_N}
        \P_N(L^{-1}b_n)b_{-n}.
    \]
    Moreover,
    \[
        \sum_{n\in\Gamma_N}
        \P_N(\Lambda^{-1}b_n)b_{-n}
        =
        \sum_{n\in\Gamma_N}\langle n\rangle^{-2}.
    \]
    Hence,
    \begin{equation}\label{eq:renormalized_variance_D}
        M_N
        =
        \sum_{n\in\Gamma_N}
        \P_N(Db_n)b_{-n}.
    \end{equation}

    For \(m\in\Z^2\), define
    \[
        \Gamma_N(m)
        \coloneqq
        \{n\in\Z^2:
            n\in\Gamma_N,\ n+m\in\Gamma_N
        \}.
    \]
    Taking the \(m\)-th Fourier coefficient in
    \eqref{eq:renormalized_variance_D} gives
    \begin{equation}\label{eq:renormalized_variance_fourier}
        \widehat{M_N}(m)
        =
        \sum_{n\in\Gamma_N(m)}
        \langle Db_n,b_{n+m}\rangle_{L^2_x}.
    \end{equation}
    The resolvent identities, understood in the quadratic-form sense,
    yield
    \begin{align*}
        D
        &=
        -L^{-1}M_\mu\Lambda^{-1}
        \\
        &=
        -\Lambda^{-1}M_\mu\Lambda^{-1}
        +
        \Lambda^{-1}M_\mu
        L^{-1}M_\mu\Lambda^{-1}.
    \end{align*}
    Therefore,
    \[
        \widehat{M_N}(m)
        =
        -a_N(m)\widehat M_\mu(m)+R_N(m),
    \]
    where
    \[
        a_N(m)
        \coloneqq
        \sum_{n\in\Gamma_N(m)}
        \langle n\rangle^{-2}
        \langle n+m\rangle^{-2}
    \]
    and
    \begin{align*}
        R_N(m)
        \coloneqq
        \sum_{n\in\Gamma_N(m)}
        &\langle n\rangle^{-2}
        \langle n+m\rangle^{-2}
        \left\langle
            L^{-1}(M_\mu b_n),
            M_\mu b_{n+m}
        \right\rangle_{L^2_x}.
    \end{align*}

    We first consider the linear term. Define
    \[
        a(m)
        \coloneqq
        \sum_{n\in\Z^2}
        \langle n\rangle^{-2}
        \langle n+m\rangle^{-2}.
    \]
    Note, the following standard convolution estimate
    \begin{align}\label{eq:conv_ineq}
        \sum_{n \in \mathbb{Z}^2} \langle n \rangle^{-\alpha} \langle n + m \rangle^{-\alpha} \lesssim 
        \begin{cases}
             \jb{m}^{2-2\alpha} \quad \text{if $\alpha \in (1,2)$}
             \\
             \jb{m}^{-2} (1+\ln(\jb{m}^2)) \quad \text{if $\alpha = 2$.}
        \end{cases}
    \end{align}
    A proof is provided in Appendix Lemma \ref{lemma:conv_ineq}.
    This convolution estimate yields
    \begin{align}\label{eq:a_m_uniform_bound}
        \sup_{m\in\Z^2}a(m)<\infty.
    \end{align}
    For every fixed \(m\), we have
    \[
        a_N(m)\xlongrightarrow[]{N \to \infty} a(m),
        \qquad \text{and} \qquad
        0\leq a_N(m)\leq a(m).
    \]
    In particular, $\sup_m\langle m\rangle^\sigma a(m)<\infty$,
    so the assumption $M\in H^{-\sigma}(\T^2)$ implies
    \[
        \sum_{m\in\Z^2}
        a(m)^2|\widehat M_\mu(m)|^2
        \lesssim
        \|M_\mu\|_{H^{-\sigma}_x}^2<\infty.
    \]
    Dominated convergence therefore gives
    $a_N\widehat M\to a\widehat M$ in $\ell^2(\Z^2)$.

    It remains to treat \(R_N\). Multiplication by a Fourier mode
    satisfies
    \begin{equation}\label{eq:modulation_Hminus_sigma}
        \|M_\mu B_k\|_{H^{-\sigma}_x}
        \lesssim
        \langle k\rangle^\sigma
        \|M_\mu\|_{H^{-\sigma}_x}.
    \end{equation}
    Indeed, this follows directly from the standard algebraic inequality $\langle x+y \rangle^s \le 2^{-\frac{s}{2}} \langle x \rangle^s \langle y \rangle^{-s}$ for $s < 0$ after shifting Fourier frequencies.

    Next, we note that we obtain sufficient smoothing through the operator $L$: By Corollary~\ref{cor:equ_norms}, we have that $\Lambda^s$ and $L^s$ are form equivalent for all $s\in [-1,1]$.
    Furthermore, since $L^{-s}$ is a bounded for $s \geq 0$ and a self-adjoint operator on $L^2(\T^2)$, we have that $\langle L^{-s} u, u \rangle_{L^2_x} \les \|u\|^2_{L^2_x}$.
    Combined with \eqref{eq:modulation_Hminus_sigma}, the duality paring of $H^{\sigma}(\T^2)$ and $H^{-\sigma}(\T^2)$, these considerations yield:
    \begin{align*}
        \left|
        \left\langle
            L^{-1}(M_\mu b_n),
            M_\mu b_{n+m}
        \right\rangle
        \right|
        &\leq
        \|L^{-1}(M_\mu b_n)\|_{H^\sigma_x}
        \|M_\mu b_{n+m}\|_{H^{-\sigma}_x}
        \\
        &\lesssim
        \|M_\mu b_n\|_{H^{-1}_x}
        \|M_\mu b_{n+m}\|_{H^{-\sigma}_x}
        \\
        &\lesssim
        \langle n\rangle^\sigma
        \langle n+m\rangle^\sigma
        \|M_\mu\|_{H^{-\sigma}_x}^2.
    \end{align*}
    Hence, using the convolution estimate in \eqref{eq:conv_ineq} with $\alpha = 2-\sigma$, we have uniformly in \(N \in \N\),
    \begin{align}\label{eq:R_N_uniform_bound}
        |R_N(m)|
        &\lesssim
        \|M_\mu\|_{H^{-\sigma}_x}^2
        \sum_{n\in\Z^2}
        \langle n\rangle^{-2+\sigma}
        \langle n+m\rangle^{-2+\sigma}
        \nonumber\\
        &\lesssim
        \|M_\mu\|_{H^{-\sigma}_x}^2
        \langle m\rangle^{-2+2\sigma}.
    \end{align}
    Thus, the full series defining $R(m)$ is absolutely
    convergent for every $m$, and $R_N(m)\to R(m)$.
    Since $\sigma<\frac12$,
    \[
        \sum_{m\in\Z^2}\langle m\rangle^{-4+4\sigma}<\infty.
    \]
    Dominated convergence consequently gives
    $R_N\to R$ in $\ell^2(\Z^2)$.

    It follows by Parseval's identity that $M_N$ converges
    in $L^2(\T^2)$ to $M_\mu$. Hence, $M_\mu\in L^2(\T^2)$ and
    $M_N\to M_\mu$ strongly in $L^2(\T^2)$.

    Finally, for
    \[
        \widetilde M_N
        =
        \sum_{n\in\Gamma_N}(Db_n)b_{-n},
    \]
    the same Fourier computation applies with
    $\Gamma_N(m)$ replaced by $\Gamma_N$ in the definitions
    of $a_N(m)$ and $R_N(m)$.
    The pointwise limits remain $a(m)$ and $R(m)$, and
    all the preceding bounds remain valid.
    The same dominated-convergence argument therefore
    gives $\widetilde M_N\to M$ in $L^2(\T^2)$, proving
    both limits in \eqref{eq:M_mu_expansion}.

\end{proof}

Finally, we are prepared to prove the conditional uniqueness of the invariant measure of the SDMFW.
\begin{proof}[Proof of Theorem \ref{theorem:uniqueness_SDMFW}]
    Suppose there exists an invariant measure $\mu \in \mathcal{W}$ of the \eqref{eq:SDMFW}, and set $M_\mu=\mathbb E_\mu[\wick{u^2}]$ and $L=\Lambda+M_\mu$. 
    Due to the equivalence of Hilbert spaces established in Corollary~\ref{cor:equ_norms} we can apply Proposition~\ref{prop_characterization_invmsr} to infer that $L$ is strictly positive and that $\mu = \mathcal{N}(0, L^{-1}) \otimes \mathcal{N}(0, \varepsilon^{-2}I)$. 
    
    Applying Lemma \ref{lem:renormalized_variance_convergence}
    to the position marginal, we obtain
    $M_\mu\in L^2(\T^2)$ and
    \[
        \widetilde M_N
        \coloneqq
        \sum_{n\in\Gamma_N}
        (L^{-1}-\Lambda^{-1})b_n\,b_{-n}
        \longrightarrow M_\mu
        \qquad\text{in }L^2(\T^2).
    \]

    Since $M_\mu$ is real-valued, the resolvent identity
    $L^{-1}-\Lambda^{-1}=-L^{-1}M_\mu\Lambda^{-1}$
    gives, for every $N$,
    \begin{align*}
        \langle\widetilde M_N,M_\mu\rangle_{L^2_x}
        &=
        \sum_{n\in\Gamma_N}
        \left\langle
            (L^{-1}-\Lambda^{-1})b_n,M_\mu b_n
        \right\rangle_{L^2_x}\\
        &=
        -\sum_{n\in\Gamma_N}\langle n\rangle^{-2}
        \left\langle
            L^{-1}(M_\mu b_n),M_\mu b_n
        \right\rangle_{L^2_x}\\
        &\leq0.
    \end{align*}
    Here all sums are finite, $M_\mu b_n\in L^2(\T^2)$,
    and the last inequality follows from positivity of
    $L^{-1}$.
    Passing to the limit using the strong $L^2$ convergence
    of $\widetilde M_N$ yields
    \[
        \|M_\mu\|_{L^2_x}^2
        =
        \lim_{N\to\infty}
        \langle\widetilde M_N,M_\mu\rangle_{L^2_x}
        \leq0.
    \]
    Consequently, $M_\mu=0$, so $L=\Lambda$ and
    $\mu=\rho_\eps$.

    Conversely, the chosen Wick ordering gives
    $\mathbb E_{\rho_\eps}[\wick{u^2}]=0$.
    The stationary linear damped-wave dynamics with
    position covariance $\Lambda^{-1}$ and velocity
    covariance $\eps^{-2}I$ therefore provide a stationary
    solution of \eqref{eq:SDMFW} with law $\rho_\eps$.
    Hence, $\rho_\eps$ is the unique invariant measure
    in $\mathcal W$.
\end{proof}

\section{Small-mass behavior of the invariant measures}\label{sec_invdynamics_SKA}

Recall that
\[
    \rho_\eps
    =
    \mathcal{N}(0, \Lambda^{-1}) \otimes \mathcal{N}(0, \varepsilon^{-2}I),
    \qquad \text{and} \qquad
    \rho_0=\mathcal N(0,\Lambda^{-1}).
\]
The uniqueness of the invariant measure $\rho_0$ of the \eqref{eq:SMFH} follows analogous to the uniqueness of the invariant measure of the \eqref{eq:SDMFW}.
Concerning the Smoluchowski-Kramers type behavior of these invariant measures, the position marginal of \(\rho_\eps\) equals \(\rho_0\) for every
\(\eps>0\). By contrast, the velocity marginal escapes to infinity as \(\eps\to0\). This will be demonstrated rigorously in the following proof.

\begin{proof}[Proof of Theorem \ref{theorem:convergence_uniqueness_SMFH}]
    We first establish uniqueness of the invariant measure of the
\eqref{eq:SMFH}. Since this proof follows the exact same steps presented in Section~\ref{Chpt_Invariant_Dynamics} we only discuss the important steps. 
Let $\mu\in\mathcal W_0$ be an invariant measure,
and set $M_\mu\coloneqq\mathbb E_\mu[\wick{u^2}]$ and
$L\coloneqq\Lambda+M_\mu$. By stationarity, the equation reduces to
\[
    du(t)=-Lu(t)\,dt+\sqrt{2}\,dW_t.
\]
Lemmas \ref{lem:relative_bound} and \ref{lemma_equ_ev} provide a
self-adjoint realization of $L$ with compact resolvent. If
$Le_j=\lambda_j e_j$ in a real orthonormal eigenbasis, the
corresponding coordinates satisfy
\[
    d\langle u(t),e_j\rangle
    =-\lambda_j\langle u(t),e_j\rangle\,dt+\sqrt{2}\,d\beta_j(t),
\]
where the $\beta_j$ are independent standard Brownian motions.
The existence of an invariant probability measure excludes
$\lambda_j\leq0$, since the corresponding scalar equation has
no invariant probability measure. Thus, $L$ is strictly positive.
Every finite collection of these coordinates consequently has
the unique invariant law given by the product of the centered
Gaussian measures with variances $\lambda_j^{-1}$. It follows
that $\mu=\mathcal N(0,L^{-1})$, with covariance understood in
the distributional sense. The remainder of the argument is
identical to the self-consistency argument in the proof of
Theorem \ref{theorem:uniqueness_SDMFW}: Lemma
\ref{lem:renormalized_variance_convergence} yields
$M_\mu\in L^2(\T^2)$, and the resolvent identity together with
positivity of $L^{-1}$ gives $M_\mu=0$. Hence, $L=\Lambda$ and
$\mu=\rho_0$. Conversely, the chosen Wick ordering gives
$\mathbb E_{\rho_0}[\wick{u^2}]=0$, so the stationary
Ornstein-Uhlenbeck process with drift $-\Lambda$ and covariance
$\Lambda^{-1}$ is a stationary solution of \eqref{eq:SMFH}.
Therefore, $\rho_0$ is the unique invariant measure in
$\mathcal W_0$.

Next, we consider the stated Smoluchowski-Kramers type behavior.
    The identity
    \[
        (\pi_1)_*\rho_\eps=\rho_0
    \]
    follows immediately from the definition of \(\rho_\eps\).

    It remains to consider the velocity marginal. Let
    \(\psi_1\) be a centered Gaussian random distribution whose
    covariance form on \(L^2(\T^2)\) is the identity. 
    Let $\{h_n\}_{n \in \Z^2}$ be a family of random variables that are independent for $n \in (\Z_+\times \{0\}) \cup (\Z\times \Z_+) \cup\{(0, 0)\}$ such that for $n \neq 0$, $h_n$ is a complex-valued centered Gaussian random variable with variance $2$ satisfying $h_{-n} = \overline{h_n}$, and for $n = 0$, $h_0$ is a real-valued centered Gaussian random variable with variance $1$.
    Using these Gaussian variables, it may be realized as
    \[
        \psi_1(x)
        =
        h_0
        +
        \frac{1}{\sqrt{2}}
        \sum_{n\in\Z^2\setminus\{0\}}
        h_ne^{in\cdot x}.
    \]
    Since \(\sigma>0\),
    \[
        \mathbb E\big[
            \|\psi_1\|_{H^{-1-\sigma}_x}^2
        \big]
        =
        \sum_{n\in\Z^2}
        \langle n\rangle^{-2-2\sigma}
        <\infty.
    \]
    Hence, the series defining \(\psi_1\) converges in
    \(L^2(\Omega;H^{-1-\sigma}(\T^2))\).

    For \(\eps>0\), define
    \[
        \psi_\eps\coloneqq\eps^{-1}\psi_1.
    \]
    Then \(\psi_\eps\) has covariance form \(\eps^{-2}I\), and
    therefore
    \[
        \operatorname{Law}(\psi_\eps)
        =
        (\pi_2)_*\rho_\eps.
    \]
    Moreover,
    \[
        \|\psi_\eps\|_{H^{-1-\sigma}_x}
        =
        \eps^{-1}
        \|\psi_1\|_{H^{-1-\sigma}_x}.
    \]
    The Gaussian field \(\psi_1\) is nonzero almost surely. Indeed,
    \(\psi_1=0\) would in particular imply \(h_0=0\), which has
    probability zero. Consequently,
    \[
        \|\psi_\eps\|_{H^{-1-\sigma}_x}
        \longrightarrow\infty
        \qquad\text{almost surely as }\eps\to0.
    \]

    Let
    \[
        f\in
        C_{\mathrm{rv}}
        (H^{-1-\sigma}(\T^2)).
    \]
    By the radial-vanishing property,
    \[
        f(\psi_\eps)\longrightarrow0
        \qquad\text{almost surely}.
    \]
    Since \(f\) is bounded, the dominated convergence theorem gives
    \begin{align*}
        \int_{H^{-1-\sigma}(\T^2)}
        f(v)\,d((\pi_2)_*\rho_\eps)(v)
        = \mathbb E[f(\psi_\eps)]
        \longrightarrow0.
    \end{align*}
    This proves the asserted convergence.
\end{proof}
\appendix
\section{Appendix}

\begin{lemma}\label{lemma:conv_ineq}
    Let $m \in \mathbb{Z}^2$ and $\alpha \in (1, 2)$. Then we have the bound
    \begin{equation*}
        S(m) \coloneqq \sum_{n \in \mathbb{Z}^2} \langle n \rangle^{-\alpha} \langle n + m \rangle^{-\alpha} \lesssim_\alpha \langle m \rangle^{2-2\alpha}.
    \end{equation*}
    If $\alpha = 2$, then we have the bound
    \begin{align*}
        S_2(m) \coloneqq \sum_{n \in \mathbb{Z}^2} \langle n \rangle^{-2} \langle n + m \rangle^{-2} \lesssim \jb{m}^{-2} (1+\ln(\jb{m}^2)).
    \end{align*}
\end{lemma}

\begin{proof}
    We cover the domain of summation into three regions: $R_1 \coloneqq \{n \in \mathbb{Z}^2 : |n| \le \frac{|m|}{2}\}$, $R_2 \coloneqq \{n \in \mathbb{Z}^2 : |n+m| \le \frac{|m|}{2}\}$, and $R_3 \coloneqq \{n \in \mathbb{Z}^2 : |n| > \frac{|m|}{2} \text{ and } |n + m| > \frac{|m|}{2}\}$. Since $\mathbb{Z}^2 \subset R_1 \cup R_2 \cup R_3$, establishing the bound on each of these three subregions is sufficient to prove the statement for the entire sum.

    The bound for $R_1$ follows by applying the reverse triangle inequality. For any $n \in R_1$, we have $|n| \le \frac{|m|}{2}$, which implies:
    \begin{equation*}
        |n+m| \ge |m| - |n| \ge |m| - \frac{|m|}{2} = \frac{|m|}{2}.
    \end{equation*}
    Because $\alpha > 0$, this implies $\langle n + m \rangle^{-\alpha} \lesssim \langle m \rangle^{-\alpha}$. Approximating the remaining sum by a polar integral and utilizing the fact that $2 - \alpha > 0$ for $\alpha \in (1, 2)$, we obtain
    \begin{align}\label{S_R1_bound}
        \begin{split}
        S_{R_1}(m) &= \sum_{n \in R_1} \langle n \rangle^{-\alpha}\langle n+m \rangle^{-\alpha} 
        \lesssim \langle m \rangle^{-\alpha} \sum_{|n| \le \frac{|m|}{2}} \langle n \rangle^{-\alpha} \\
        &\lesssim \langle m \rangle^{-\alpha} \left(1+\int_0^{\frac{|m|}{2}} (1 + r^2)^{-\frac{\alpha}{2}} r \, dr \right)
        = \langle m \rangle^{-\alpha} \left(1 + \frac{1}{2 - \alpha} \left[ (1+r^2)^{1-\frac{\alpha}{2}} \right]_{0}^{\frac{|m|}{2}} \right) \\
&\lesssim_\alpha \langle m \rangle^{-\alpha}
    \left(1+\langle m \rangle^{2-\alpha}\right)
\lesssim_\alpha \langle m \rangle^{2-2\alpha}.        
\end{split}
    \end{align}
    The bound in region $R_2$ follows by shifting the summation index to $n' = n+m$. 

    To show the bound in the final region $R_3$, we apply Young's product inequality
    \begin{equation*}
        S_{R_3}(m) \le \frac{1}{2} \sum_{n \in R_3}\langle n \rangle^{-2\alpha} + \frac{1}{2} \sum_{n \in R_3}\langle n+m \rangle^{-2\alpha}.
    \end{equation*}
    By symmetry, it suffices to bound the first term. Comparing the sum to a polar integral and noting that $1 - \alpha < 0$, we obtain
    \begin{align*}
        \sum_{n \in R_3} \langle n \rangle^{-2\alpha} 
        \le \sum_{|n| > \frac{|m|}{2}} \langle n \rangle^{-2\alpha} 
        \lesssim \int_{\frac{|m|}{2}}^\infty (1+r^2)^{-\alpha} r \, dr 
        = \frac{1}{2 - 2\alpha} \left[ (1+r^2)^{1-\alpha} \right]_{\frac{|m|}{2}}^\infty 
        \lesssim \langle m \rangle^{2-2\alpha}.
    \end{align*}
    The proof for the bound of $S_2$ follows analogously, but the estimate in \eqref{S_R1_bound} changes to 
    \begin{align*}
        S_{2_{R_1}}(m) &= \sum_{n \in R_1} \langle n \rangle^{-2}\langle n+m \rangle^{-2} 
        \\
        &\lesssim \langle m \rangle^{-2} \left(1+ \int_0^{\frac{|m|}{2}} (1 + r^2)^{-1} r \, dr \right)
        \\
&\leq \langle m \rangle^{-2}
    \left(1+\ln(1+|m|^2)\right).
    \end{align*}
    This yields the desired estimate.
\end{proof}

\section*{Acknowledgements}
The author is grateful to Herbert Koch and Ruoyuan Liu for their supervision of the author's master's thesis and for introducing the author to the study of stochastic partial differential equations in infinite dimensions. Particular thanks are due to Herbert Koch for pointing out that an infinite-dimensional Sobolev space, equipped with its norm topology, admits no nonzero continuous functions with compact support. This observation informed the choice of test functions used to describe the behavior of the velocity marginals.

\section*{Disclosure of AI assistance}
During the preparation of this manuscript, the author used
OpenAI's ChatGPT for assistance with writing and exposition
and with the development and revision of parts of the
mathematical arguments, including suggestions for auxiliary
lemmas and proofs. Specific uses are documented in the
annotated LaTeX source files available in the
\href{https://github.com/alexanderczernik/alexanderczernik.github.io/blob/dee787a89a7c0ef80735113e041fa398cd73a694/assets/files/uniqueness_manuscript}
{accompanying source repository}.
The author takes full responsibility for the mathematical
statements, proofs, references, and final text.

\printbibliography

@book{DaPrato_Zabczyk_2014,
  title={Stochastic Equations in Infinite Dimensions},
  author={Da Prato, Giuseppe and Zabczyk, Jerzy},
  series={Encyclopedia of Mathematics and its Applications},
  volume={45},
  year={1992},
  publisher={Cambridge University Press},
  address={Cambridge}
}

@book {schmudgen_operators,
    AUTHOR = {Schm\"udgen, Konrad},
     TITLE = {Unbounded self-adjoint operators on {H}ilbert space},
    SERIES = {Graduate Texts in Mathematics},
    VOLUME = {265},
 PUBLISHER = {Springer, Dordrecht},
      YEAR = {2012},
     PAGES = {xx+432},
      ISBN = {978-94-007-4752-4},
   MRCLASS = {47-01 (47B25 47E05)},
  MRNUMBER = {2953553},
MRREVIEWER = {G.\ V.\ Rozenblum},
       DOI = {10.1007/978-94-007-4753-1},
       URL = {https://doi.org/10.1007/978-94-007-4753-1},
}

@misc{liu2025,
      title={Hyperbolic $O (N)$ linear sigma model and its mean-field limit}, 
      author={Ruoyuan Liu and Shao Liu and Tadahiro Oh},
      year={2025},
      eprint={2511.21950},
      archivePrefix={arXiv},
      primaryClass={math.AP},
      url={https://arxiv.org/abs/2511.21950}, 
}

@article {Shen2022,
    AUTHOR = {Shen, Hao and Smith, Scott A. and Zhu, Rongchan and Zhu,
              Xiangchan},
     TITLE = {Large {$N$} limit of the {$O(N)$} linear sigma model via
              stochastic quantization},
   JOURNAL = {Ann. Probab.},
  FJOURNAL = {The Annals of Probability},
    VOLUME = {50},
      YEAR = {2022},
    NUMBER = {1},
     PAGES = {131--202},
      ISSN = {0091-1798,2168-894X},
   MRCLASS = {60H15 (35R60 60F17 81S20)},
  MRNUMBER = {4385125},
MRREVIEWER = {Francesco\ Grotto},
       DOI = {10.1214/21-aop1531},
       URL = {https://doi.org/10.1214/21-aop1531},
}

@misc{czernik2026,
      title={On the Smoluchowski-Kramers approximation for the hyperbolic $O(N)$ linear sigma model and its mean-field limit}, 
      author={Alexander Czernik and Ruoyuan Liu and Shao Liu},
      year={2026},
      eprint={2606.15214},
      archivePrefix={arXiv},
      primaryClass={math.AP},
      url={https://arxiv.org/abs/2606.15214}, 
}

@misc{tolomeo2023ergodicityhyperbolicpphi2model,
      title={Ergodicity for the hyperbolic $P(\Phi)_2$-model}, 
      author={Leonardo Tolomeo},
      year={2023},
      eprint={2310.02190},
      archivePrefix={arXiv},
      primaryClass={math.PR},
      url={https://arxiv.org/abs/2310.02190}, 
}

@article {Parisi1981,
    AUTHOR = {Parisi, G. and Wu, Yong Shi},
     TITLE = {Perturbation theory without gauge fixing},
   JOURNAL = {Sci. Sinica},
  FJOURNAL = {Scientia Sinica. Zhongguo Kexue},
    VOLUME = {24},
      YEAR = {1981},
    NUMBER = {4},
     PAGES = {483--496},
      ISSN = {0582-236x},
   MRCLASS = {81E99},
  MRNUMBER = {626795},
}

@book{GlimmJaffe1987,
  author    = {Glimm, James and Jaffe, Arthur},
  title     = {Quantum Physics: A Functional Integral Point of View},
  edition   = {2},
  publisher = {Springer-Verlag},
  address   = {New York},
  year      = {1987},
  doi       = {10.1007/978-1-4612-4728-9}
}

@article{DaPratoDebussche2003,
  author  = {Da Prato, Giuseppe and Debussche, Arnaud},
  title   = {Strong solutions to the stochastic quantization equations},
  journal = {The Annals of Probability},
  volume  = {31},
  number  = {4},
  pages   = {1900--1916},
  year    = {2003},
  doi     = {10.1214/aop/1068646370}
}

@article{MourratWeber2017Plane,
  author  = {Mourrat, Jean-Christophe and Weber, Hendrik},
  title   = {Global well-posedness of the dynamic {$\Phi^4$}
             model in the plane},
  journal = {The Annals of Probability},
  volume  = {45},
  number  = {4},
  pages   = {2398--2476},
  year    = {2017},
  doi     = {10.1214/16-AOP1116}
}

@article{HairerMattingly2018,
  author  = {Hairer, Martin and Mattingly, Jonathan},
  title   = {The strong {Feller} property for singular stochastic {PDEs}},
  journal = {Annales de l'Institut Henri Poincar{\'e},
             Probabilit{\'e}s et Statistiques},
  volume  = {54},
  number  = {3},
  pages   = {1314--1340},
  year    = {2018},
  doi     = {10.1214/17-AIHP840}
}

@article{TsatsoulisWeber2018,
  author  = {Tsatsoulis, Pavlos and Weber, Hendrik},
  title   = {Spectral gap for the stochastic quantization equation
             on the 2-dimensional torus},
  journal = {Annales de l'Institut Henri Poincar{\'e},
             Probabilit{\'e}s et Statistiques},
  volume  = {54},
  number  = {3},
  pages   = {1204--1249},
  year    = {2018},
  doi     = {10.1214/17-AIHP837}
}

@article{DuaneEtAl1987,
  author  = {Duane, Simon and Kennedy, A. D. and
             Pendleton, Brian J. and Roweth, Duncan},
  title   = {Hybrid {Monte Carlo}},
  journal = {Physics Letters B},
  volume  = {195},
  number  = {2},
  pages   = {216--222},
  year    = {1987},
  doi     = {10.1016/0370-2693(87)91197-X}
}

@article{GKOT2022,
  author  = {Gubinelli, Massimiliano and Koch, Herbert and
             Oh, Tadahiro and Tolomeo, Leonardo},
  title   = {Global Dynamics for the Two-dimensional Stochastic
             Nonlinear Wave Equations},
  journal = {International Mathematics Research Notices},
  volume  = {2022},
  number  = {21},
  pages   = {16954--16999},
  year    = {2022},
  doi     = {10.1093/imrn/rnab084}
}

@article{OTWZ2026,
  author  = {Oh, Tadahiro and Tolomeo, Leonardo and
             Wang, Yuzhao and Zheng, Guangqu},
  title   = {Hyperbolic {$P(\Phi)_2$}-model on the Plane},
  journal = {Communications in Mathematical Physics},
  volume  = {407},
  number  = {2},
  pages   = {34},
  year    = {2026},
  doi     = {10.1007/s00220-025-05486-0}
}

@misc{Tolomeo2023Ergodicity,
  author        = {Tolomeo, Leonardo},
  title         = {Ergodicity for the hyperbolic {$P(\Phi)_2$}-model},
  year          = {2023},
  eprint        = {2310.02190},
  archivePrefix = {arXiv},
  primaryClass  = {math.PR},
  doi           = {10.48550/arXiv.2310.02190}
}

@article{Wilson1973,
  author  = {Wilson, Kenneth G.},
  title   = {Quantum field-theory models in less than 4 dimensions},
  journal = {Physical Review D},
  volume  = {7},
  year    = {1973},
  pages   = {2911--2926}
}

@misc{DelgadinoSmith2026,
  author        = {Delgadino, Mat{\'i}as G. and Smith, Scott A.},
  title         = {Mass generation for the two dimensional {$O(N)$}
                   Linear Sigma Model in the large {$N$} limit},
  year          = {2026},
  eprint        = {2601.19630},
  archivePrefix = {arXiv},
  primaryClass  = {math.PR},
  url           = {https://arxiv.org/abs/2601.19630}
}

@article{Zine2025,
  author  = {Zine, Younes},
  title   = {{Smoluchowski--Kramers} approximation for singular
             stochastic wave equations in two dimensions},
  journal = {Electronic Journal of Probability},
  volume  = {30},
  year    = {2025},
  pages   = {Paper No. 88, 49},
  doi     = {10.1214/25-EJP1337},
  url     = {https://doi.org/10.1214/25-EJP1337}
}

@article{Hairer2006,
 ISSN = {0003486X},
 URL = {http://www.jstor.org/stable/20160014},
 author = {Martin Hairer and Jonathan C. Mattingly},
 journal = {Annals of Mathematics},
 number = {3},
 pages = {993--1032},
 publisher = {Annals of Mathematics},
 title = {Ergodicity of the 2D Navier-Stokes Equations with Degenerate Stochastic Forcing},
 urldate = {2026-09-21},
 volume = {164},
 year = {2006}
}

@book{Da_Prato_Zabczyk_1996, place={Cambridge}, series={London Mathematical Society Lecture Note Series}, title={Ergodicity for Infinite Dimensional Systems}, publisher={Cambridge University Press}, author={Da Prato, G. and Zabczyk, J.}, year={1996}, collection={London Mathematical Society Lecture Note Series}}

@article{Hairer2011,
author={Hairer, M.
and Mattingly, J. C.
and Scheutzow, M.},
title={Asymptotic coupling and a general form of Harris' theorem with applications to stochastic delay equations},
journal={Probability Theory and Related Fields},
year={2011},
month={Feb},
day={01},
volume={149},
number={1},
pages={223-259},
issn={1432-2064},
doi={10.1007/s00440-009-0250-6},
url={https://doi.org/10.1007/s00440-009-0250-6}
}

@article{RehmeierRockner2025,
  author  = {Rehmeier, Marco and R{\"o}ckner, Michael},
  title   = {On Nonlinear {Markov} Processes in the Sense of {McKean}},
  journal = {Journal of Theoretical Probability},
  volume  = {38},
  year    = {2025},
  eid     = {60},
  note    = {Article 60},
  doi     = {10.1007/s10959-025-01428-7},
  url     = {https://doi.org/10.1007/s10959-025-01428-7}
}

@article {Benyi2025,
    AUTHOR = {B\'enyi, \'Arp\'ad and Oh, Tadahiro and Zhao, Tengfei},
     TITLE = {Fractional {L}eibniz rule on the torus},
   JOURNAL = {Proc. Amer. Math. Soc.},
  FJOURNAL = {Proceedings of the American Mathematical Society},
    VOLUME = {153},
      YEAR = {2025},
    NUMBER = {1},
     PAGES = {207--221},
      ISSN = {0002-9939,1088-6826},
   MRCLASS = {42B15 (42B25 46E35)},
  MRNUMBER = {4840270},
       DOI = {10.1090/proc/17007},
       URL = {https://doi.org/10.1090/proc/17007},
}

@book{Brezis2011,
  title={Functional analysis, Sobolev spaces and partial differential equations},
  author={Brezis, Haim},
  volume={2},
  number={3},
  year={2011},
  publisher={Springer}
}

\end{document}